\documentclass[12pt,reqno]{amsart}
\usepackage{graphicx}
\usepackage{amssymb}
\makeatletter
\def\proof{\@ifnextchar[{\@oproof}{\@nproof}}
\def\@oproof[#1][#2]{\trivlist\item[\hskip\labelsep
\textit{#2 Proof of\ #1.}~]\ignorespaces}
\def\@nproof{\trivlist\item[\hskip\labelsep\textit{Proof.}~]\ignorespaces}

\makeatother

\setbox0=\hbox{$+$}
\newdimen\plusheight
\plusheight=\ht0
\def\+{\;\lower\plusheight\hbox{$+$}\;}

\setbox0=\hbox{$-$}
\newdimen\minusheight
\minusheight=\ht0
\def\-{\;\lower\minusheight\hbox{$-$}\;}

\setbox0=\hbox{$\cdots$}
\newdimen\cdotsheight
\cdotsheight=\plusheight
\def\cds{\lower\cdotsheight\hbox{$\cdots$}}

\setbox0=\hbox{$+$}

\renewcommand{\(}{\left\(}
\renewcommand{\)}{\right\)}
\renewcommand{\[}{\left[}

\numberwithin{equation}{section}
 \theoremstyle{plain}
\newtheorem{theorem}{Theorem}[section]

\begin{document}
\title[The Pythagoras' Theorem] {The Pythagoras' Theorem}
\author{Manjil P.~Saikia}
\address{Department of Mathematical Sciences, Tezpur University, Napaam, Sonitpur, Pin-784028, India}
\email{manjil@gonitsora.com}


\vspace*{0.5in}
\begin{center}
{\bf The Pythagoras' Theorem}\\[5mm]
{\footnotesize  MANJIL P.~SAIKIA}\\[3mm]


\end{center}
\vskip 5mm \noindent{\footnotesize{\bf Abstract.} We give a brief historical overview of the famous Pythagoras' theorem and Pythagoras. We present a simple proof of the result and dicsuss one direction of extension which has resulted in a famous result in number theory. We follow \cite{thales}, \cite{wiki} and \cite{wiki2} for the historical comments and sources.}

\vskip 3mm


\vskip 3mm

\noindent {\footnotesize 2010 Mathematical Reviews Classification
Numbers: 01A17, 01A20, 01A32.}

\section{{Introduction}}

In mathematics, the \textbf{\textit{Pythagoras' theorem}} is a famous result in triangle geometry. Although there are many variants of this result in various branches of mathematics, we shall restrict ourselves to the following simple version that was originally given by Pythagoras.

\begin{theorem}[Pythagoras]\label{pt}
The square of the hypotenuse (the side opposite the right angle) is equal to the sum of the squares of the other two sides.
\end{theorem}

Thus if $a$ denotes the hypotenuse of the right angled triangle and $b$ and $c$ denotes the other two sides, then the theroem says that $a^2=b^2+c^2$. Triples like $(a,b,c)$ which satisfy this theorem are known as \textit{Pythagorean triples}. It is an easy exercise to show that there are infinitely many such triples of numbers. Indeed if $(a,b,c)$ is an example, then $(ka,kb,kc)$ is also an example for any non-negative integer $k$.

Before proceeding to give a proof of Theorem \ref{pt}, we shall give a brief historical overview of the theorem and Pythagoras.

This theorem is named after the Greek mathematician \textbf{\textit{Pythagoras (ca. 570 BC -- ca. 495 BC)}}, who by tradition is credited with its proof, although it is often argued that knowledge of the theorem predates him. There is evidence that Babylonian mathematicians understood the formula, although there is little surviving evidence that they used it in a mathematical framework.

Pythagoras of Samos is often described as the first pure mathematician. He was an extremely important figure in the development of mathematics yet we know relatively little about his mathematical achievements. Unlike many later Greek mathematicians, where we have at least some of the books which they wrote, we have nothing of Pythagoras's writings. The society which he led, half religious and half scientific, followed a code of secrecy which certainly means that today Pythagoras is a mysterious figure. Whatever we know about Pythagoras is through the writings of other Greek contemporaries and later philosophers and mathematicians. Pythagoras' is considered to be the founder of a group called the Pythagorean brotherhood. In fact it is believed that Pythagoras was the first person to use the word \textit{mathematics}.
\begin{center}
  \includegraphics[height=2.5in]{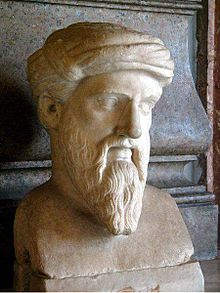}
\end{center}

\begin{center}
   \textbf{Pythagoras of Samos}
\end{center}

Pythagoras apart from being an accomplished philosopher and mathematician also dealt with music and notes. Pythagoras is believed to be the first person who systematically studied the musical notes and hence is said to be the father of music. Many of the results and teachings that Pythagoras is attributed with may have been discovered by his pupils and followers. It is believed that whatever the \textit{Pythagorean brotherhood} or the \textit{Pythagoreans} discovered, they gave credit to their master. Thus, it is matter of debate whether we must really allow so many discoveries to be credited with Pythagoras. Nonetheless, he stands as an important figure in the history of mathematics and his works and teachings have had an enormous amount of influence on the way Greek mathematics developed.

\section{{History of Theorem \ref{pt}}}

The most famous result in mathematics is perhaps the Pythagoras theorem. Every high school student if asked to state one mathematical result correctly, would invariably choose this theorem. However there is a considerable debate whether the Pythagorean theorem was discovered once, or many times in many places.

The history of the theorem can be divided into four parts: knowledge of Pythagorean triples, knowledge of the relationship among the sides of a right triangle, knowledge of the relationships among adjacent angles, and proofs of the theorem within some deductive system.

Bartel Leendert van der Waerden (1903 -- 1996) conjectured that Pythagorean triples were discovered algebraically by the Babylonians. Written between 2000 and 1786 BC, the Middle Kingdom Egyptian papyrus Berlin 6619 includes a problem whose solution is the Pythagorean triple 6:8:10, but the problem does not mention a triangle. The Mesopotamian tablet Plimpton 322, written between 1790 and 1750 BC during the reign of Hammurabi the Great, contains many entries closely related to Pythagorean triples.

In India, the Baudhayana Sulba Sutra, the dates of which are given variously as between the 8th century BC and the 2nd century BC, contains a list of Pythagorean triples discovered algebraically, a statement of the Pythagorean theorem, and a geometrical proof of the Pythagorean theorem for an isosceles right triangle. The Apastamba Sulba Sutra (ca. 600 BC) contains a numerical proof of the general Pythagorean theorem, using an area computation. Van der Waerden believed that "it was certainly based on earlier traditions". Boyer (1991) thinks the elements found in the Śulba-sũtram may be of Mesopotamian derivation.

With contents known much earlier, but in surviving texts dating from roughly the first century BC, the Chinese text Zhou Bi Suan Jing, (The Arithmetical Classic of the Gnomon and the Circular Paths of Heaven) gives a reasoning for the Pythagorean theorem for the (3, 4, 5) triangle—in China it is called the \textit{Gougu Theorem}. During the Han Dynasty (202 BC -- 220 AD), Pythagorean triples appear in The Nine Chapters on the Mathematical Art, together with a mention of right triangles. Some believe the theorem arose first in China, where it is alternatively known as the \textit{Shang Gao Theorem}, named after the Duke of Zhou's astronomer and mathematician, whose reasoning composed most of what was in the Zhou Bi Suan Jing.

Pythagoras, whose dates are commonly given as 570 -- 495 BC, used algebraic methods to construct Pythagorean triples, according to Proclus's commentary on Euclid. Proclus, however, wrote between 410 and 485 AD. According to Sir Thomas L. Heath (1861 -- 1940), no specific attribution of the theorem to Pythagoras exists in the surviving Greek literature from the five centuries after Pythagoras lived. However, when authors such as Plutarch and Cicero attributed the theorem to Pythagoras, they did so in a way which suggests that the attribution was widely known and undoubted. Whether this formula is rightly attributed to Pythagoras personally or not, but one can safely assume that it belongs to the very oldest period of Pythagorean mathematics. In fact, it is believed that Pythagoras travelled widely, and he may have actually learned about the theorem from the Babylonians.

Around 400 BC, according to Proclus, Plato gave a method for finding Pythagorean triples that combined algebra and geometry. Around 300 BC, in Euclid's Elements, the oldest extant axiomatic proof of the theorem is presented.

\section{{Proof of Theorem \ref{pt}}}

The Pythagoras' theorem is perhaps the most proved mathematical result with hundreds of known proofs. \cite{knot} contains a collection of almost 100 different proofs. The oldest dating back to the days of Euclid and the newest from the 21st century. The proofs are all varied, some of them are geometrical, some of them are algebraic, and infact there are some which uses principles of physics and differential calculus to prove this result. Below we present one way to prove the theorem.

\begin{proof}

We take four copies of a right triangle with sides a, b and c, arranged inside a square with side c as shown in the figure below. This will leave a square with side $(b-a)$ in the middle empty, while all the triangles will be similar with aread $\frac{1}{2}ab$.

\begin{center}
  \includegraphics[height=2.5in]{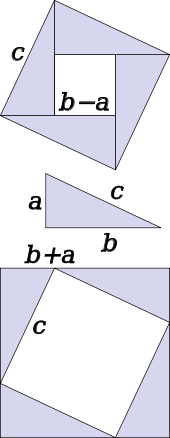}
\end{center}

The area of the large square is therefore $(b-a)^2+\frac{4}{2}ab=a^2+b^2$. This in turn is equal to $c^2$, since the square is of side length $c$.
\end{proof}

\section{{Extensions of Theorem \ref{pt}}}

A good theorem is the one which inspires other good results. In that sense, the Pythagoras' theorem has been a precursor of many wonderful mathematical ideas. The Pythagoras' Theorem can be extended to many different areas of mathematics, including but not limited to inner product spaces, non-Euclidean geometry, trigonometry, etc. Here we just mention one such aspect of how the Pythagoras' theorem has given rise to mathematics in other areas.

A Pythagorean triple $(a,b,c)$ is called \textit{primitive} if $\gcd{(a,b,c)}=1$. The following fact about primitive Pythagorean triple was found in the clay tablet Plimpton 322 (1900 -- 1600 BC) of the Bablylonians.

\begin{theorem}
  Suppose $u$ and $v$ are relatively prime positive integers, assume that not both are odd and that $u>v$. Then, if $a=2uv$, $b=u^2-v^2$ and $c=u^2+v^2$, we have $\gcd{(a,b,c)}=1$ and $a^2+b^2=c^2$.
\end{theorem}

In fact, the converse of this result is also true. The earliest record of this work is found in the famous book of \textit{Diophantus}, Arithmetica (250 AD). This brings us to one of the most famous result in number theory called the Fermat's Last Theorem. A standard restatement is,

\begin{theorem}[Fermat's Last Theorem] There are no solutions to the following problem with $(X,Y,Z)$ integers

$$X^p+Y^p+Z^p=0$$

where $XYZ\neq0$ and $p\geq3$ is a prime.

\end{theorem}

Sometime in the middle of the $17^{th}$ century Pierre de Fermat (1601 -- 1655), an amateur French mathematician wrote in the margin of his copy of Bachet's edition of the works of Diophantus,
\begin{quote}It is impossible to separate a cube into two cubes, or a biquadrate
into two biquadrates, or in general any power higher than the second
into powers of like degree; I have discovered a truly marvelous proof,
which this margin is too small to contain.\end{quote}
This is the celebrated Fermat's Last Theorem, which in modern language translates to,
\begin{theorem}\label{flt}
If $n$ is any natural number greater than $2$, the equation $$X^n+Y^n=Z^n$$ has no solutions in integers, all different from $0$.\end{theorem}

Whatever marvelous proof Fermat had for his theorem nobody found out because in all of Fermat's letters to other mathematicians he never mentioned it. His theorem gathered much publicity after his death and remained the most prized unsolved problem of mathematics for more than three centuries. Euler, the legendary Swiss mathematician proved the theorem for $n=3$, Fermat himself solved the problem for $n=4$ by using his famous method of \textit{infinite descent} and this was followed by Sophie Germain with her proof of the theorem for relatively small primes. The case when $n=5$ were proved with the efforts of Dirichlet and Legendre. Dirichlet further proved the theorem for $n=14$. After this the mathematician Lam\'{e} disposed off with the case $n=7$. This was then followed by Kummer's marvellous achievement, where he proved the result for many values. A brief survey of this result can be found in \cite{mps}. This result was finally proved in 1994 by Andrew Wiles, following the works of many mathematicians like K. A. Ribet, G. Y. Taniyama, G. Shimura, B. Mazur, G. Frey and R. Taylor.

Thus, we see how a simple result like Theorem \ref{pt} has given rise to so much beautiful mathematics which has culminated in Theorem \ref{flt}. This is really the hallmark of a good result.

\section{{Acknowledgements}}

The author wishes to thank Sirat Sandil for a careful reading of the manuscript and pointing out many grammatical and typographical errors.


\begin{thebibliography}{99}
\bibitem{thales}{W. S. Anglin and J. Lambek}, \emph{The Heritage of Thales}, UTM, Springer-Verlag, New York, 1995.
\bibitem{knot}{Cut the Knot}, \emph{http://www.cut-the-knot.org/}, accessed on 28th September, 2013.
\bibitem{mps}{M. P. Saikia}, \emph{A study of Kummer's proof of the Fermat's Last Theorem for Regular Primes}, IISER Mohali Report, 2011.
\bibitem{wiki}{Wikipedia}, \emph{Entry on Pythagoras}, accessed on 28th September, 2013.
\bibitem{wiki2}{Wikipedia}, \emph{Entry on Pythagoras' Theorem}, accessed on 28th September, 2013.
\end{thebibliography}
\end{document}